\documentclass[12pt,twoside]{amsart}
  \usepackage[latin1]{inputenc}
  \usepackage[T1]{fontenc}
  \usepackage{amsmath,amscd,amssymb,latexsym, a4}
  \usepackage[english,frenchb]{babel}
  \usepackage[all]{xy}
\usepackage{amsthm}

\newtheorem*{theorem}{Theorem}

\newtheorem*{e-proposition}{Proposition}

\newtheorem*{e-definition}{Definition\rm}

\newtheorem{theoreme}{Th\'eor\`eme}

\newtheorem{proposition}[theoreme]{Proposition}

\newtheorem*{exemple}{\it Exemple\/}

\def\og{\leavevmode\raise.3ex\hbox{$\scriptscriptstyle\langle\!\langle$~}}
\def\fg{\leavevmode\raise.3ex\hbox{~$\!\scriptscriptstyle\,\rangle\!\rangle$}}

\def\La{\Lambda }

\def \H{{\rm H}}
\def \Z{{\mathbb Z}}

\def\O{{\Omega}}

\begin{document}
\author{Vincent Franjou}
\thanks{UMR 6629 Université de Nantes/CNRS}
\email{franjou@math.univ-nantes.fr}
\address{Laboratoire Jean-Leray, Universit\'e de Nantes\\
2, rue de la Houssini\`ere\\ BP 92208\\
44322 Nantes cedex 3, France}
\title{Cohomologie de de {R}ham entière}
\maketitle{}
\selectlanguage{francais}
\begin{abstract}
On d\'ecrit la suite spectrale de Bockstein issue du complexe de
de Rham sur les entiers. L'isomorphisme de Cartier intervient
comme un endomorphisme du complexe de de Rham modulo $p$ qui
identifie les pages successives de la suite spectrale de
Bockstein. On en d\'eduit la cohomologie de de Rham enti\`ere des
espaces affines.
\end{abstract}
\bigskip\bigskip
\selectlanguage{francais}

 Pour un
anneau commutatif $A$, on note $\O_A^1$ le module des formes
différentielles absolues de degré 1 de $A$, et on note
$\O_A^*=\La_A^*(\O_A^1)$
son algèbre extérieure sur $A$. Dans l'anneau
gradué $\O_A^*$, on définit de la manière usuelle une dérivation
$d$ de degré $1$ et de carré nul. On obtient ainsi le complexe de
de Rham de $A$. Cette note étudie le cas élémentaire où l'anneau
$A$ est un anneau de polynômes à coefficients entiers, la
diff\'erentielle $d$ \'etant alors donn\'ee par la d\'erivation des
polyn\^omes. Afin de
mettre en exergue les propriétés de naturalité, on considère l'anneau
$A$
comme l'algèbre symétrique sur un groupe abélien libre de rang
fini $L$, graduée par le degré polynômial~:
$$A=S^*(L)=\bigoplus_{d\geq 0} S^d(L)$$
($S^d(L)$ est form\'e des orbites de $L^{\otimes d}$ sous
l'action de permutation des $d$ facteurs).
D\'esignant par $\La^i(L)$ la $i$-ième puissance extérieure du
groupe ab\'elien $L$, on a~: $\O_A^i=A\otimes\La^iL$. Posant :
$$(\O_A^i)_n=S^{n-i}(L)\otimes\La^i(L)=:\O^i_n(L)\ ,$$
on voit que le complexe de de Rham est lui-même gradué. Les
propriétés de naturalité de sa différentielle $d$ permettent de le
considérer comme un complexe de foncteurs de $L$. On oublie donc
la mention du groupe abélien libre $L$, et la notation
$\O_n^{\star}$ désigne dans la suite le complexe de foncteurs
donné par le complexe de de Rham en degré total $n$ :
  $$\O_n^{\star}\ :\ \xymatrix{S^n\ar[r]^-{d} &
  S^{n-1}\otimes\La^1\ar[r]^-{d}&\cdots\
  S^{n-i}\otimes\La^i\ar[r]^-{d}&\cdots\ar[r]^-{d}& \La^n\ar[r]^-{d}&
  0\ar[r]^-{d}&0 \ \cdots}\ .$$
On note $\H^i(\O_n^{\star})$ le groupe de cohomologie
correspondant (c'est un foncteur de $L$).\par Le fait que la
cohomologie de de Rham de l'espace affine ainsi d\'efinie soit non nulle est un
phénomène bien connu, étudié dans une note de Pierre Cartier
\cite{Cartier}. Ce ph\'enom\`ene est à la base de travaux sur la suite
spectrale de Hodge vers de Rham \cite{Deligne}. Il permet aussi
d'effectuer efficacement des calculs de
cohomologie des foncteurs \cite{FLS}. Ces travaux omettent pourtant de
signaler le calcul de la cohomologie de de Rham dans le cas des
entiers, ce qui donne prétexte à cette note.
\bigskip
\begin{proposition}\label{torsion}
La cohomologie de de Rham de degré total $n$ est un groupe fini
annulé par $n$.
\end{proposition}
\par\noindent\emph{D\'emonstration.} L'algèbre graduée $\O^*$ est
aussi muni d'une différentielle de Koszul $\kappa$, qui est une
dérivation naturelle envoyant les polynômes sur $0$, et une forme
différentielle de degré $1$, $dx$, sur le polynôme $x$. Elle est
liée à la différentielle de de Rham $d$ par la formule d'Euler :
$d\kappa +\kappa d=n$. Le résultat en découle.
\par\medskip
Puisque la cohomologie de de Rham est un groupe de torsion, elle est
d\'etermim\'ee, en tant que groupe ab\'elien, par la suite
spectrale de Bockstein du complexe de de Rham pour chaque nombre premier $p$.
Soit $p$ un nombre premier. Pour obtenir la suite spectrale de Bockstein, on
considère la suite exacte courte de complexes :
$$0\to\O_{n}^{\star}\to\O_{n}^{\star}\to \O_n^{\star}\otimes\Z/p\Z\to 0\ ,$$
induite par la multiplication par $p$. Elle engendre un couple
exact
$$  \xymatrix {
    \H (\O _n^{\star}) \ar[rr]^-p  && \H (\O_n^{\star})\ar[ld] \\
    & \H (\O_{n}^{\star}\otimes\Z/p)\ar[lu]^\partial }$$
dont les couples dérivés fournissent les pages successives de la
suite spectrale de Bockstein. \par Choisissons une base de
$\O^1_1$. Le remplacement des variables par leur puissance
$p$-ième permet de définir un morphisme de complexes $F$ (non
naturel), qui envoie $x$ sur $x^p$ et $dx$ sur $px^{p-1}dx$. On
l'appellera Frobenius.
\begin{proposition}
Soit $p$ un nombre premier. Le Frobenius induit une application
naturelle :
$$
F_*\ :\ \H^i(\O_n^{\star})\to \H^i(\O_{pn}^{\star})\ .
$$
\end{proposition}
\par\medskip
On définit aussi un endomorphisme de l'algèbre $\O$ (non naturel)
qui envoie $x$ sur $x^p$ et $dx$ sur $x^{p-1}dx$. Il définit un
morphisme d'algèbres $C^{-1}$, naturel, de $\O\otimes\Z/p$ vers
$\H(\O^{\star}\otimes\Z/p)$, qui est un isomorphisme
\cite{Cartier}. C'est son isomorphisme réciproque que l'on appelle
isomorphisme de Cartier.
\begin{theoreme}
Soit $p$ un nombre premier. La suite spectrale de Bockstein du
complexe de de Rham est stationnaire, chaque page s'identifiant,
par l'isomorphisme de Cartier, au complexe de de Rham sur le corps
à $p$ éléments.
\end{theoreme}
\par\noindent\emph{D\'emonstration.} La différentielle $d_1$ de la
première page s'obtient en composant le connectant $\partial$ et
la reduction mod. $p$. Pour l'identifier, un calcul en bas degré
suffit. Alternativement, on peut s'aider d'une base de $\O^1_1$
pour effectuer le calcul suivant. Considérons une forme $\omega =
Pdx_1\dots dx_i$ dans $\O^i$, et notons la réduction mod. $p$ par
une barre. Par définition, $C^{-1}(\overline{\omega})$ est la
classe de cohomologie de la réduction mod. $p$ de
$\frac{F(\omega)}{p^i}$. On a donc :
$$\partial C^{-1}(\overline{\omega})=
\partial[\overline{\frac{F(\omega)}{p^i}}]=[\frac{dF(\omega)}{p^{i+1}}]$$
dont la réduction mod. $p$ n'est autre que
$C^{-1}(d\overline{\omega})$. La première page de la suite
spectrale de Bockstein est donc le complexe de de Rham modulo
$p$.\par
Pour les autres pages, on a le résultat suivant~:
\begin{proposition} Le Frobenius induit un morphisme du couple exact de Bockstein
sur son couple dérivé :
$$  \xymatrix {
    \H (\O^{\star}) \ar[rr]^-p  \ar[dd]^{F_*} && \H (\O^{\star}) \ar[dd]^{F_*}\ar[ld] \\
    & E_1\ar[lu]^\partial \ar[dd] \\
    p\H (\O^{\star}) \ar[rr] && p\H (\O^{\star})\ar[ld] \\
    & E_2 \ar[lu]^\partial }\ .$$
Le morphisme $\H (\O^{\star} /p\O^{\star})=E_1\to E_2$ est un
isomorphisme induit par l'isomorphisme $C^{-1}$ :
$\O_n\otimes\Z/p\Z\to\H(\O_{pn}^{\star}\otimes\Z/p\Z)$.
\end{proposition}
\par\noindent
La seule chose qui mérite une vérification est que l'image de
$F_*$ est divisible par $p$. Prenons donc un cocycle $\omega$ dans
$\O_n$. Comme $n\omega$ est un cobord (proposition \ref{torsion}),
c'est aussi le cas de $F(n\omega)$. Il en résulte que
$\displaystyle\frac{F(n\omega)}{np}$ est un cocycle de $\O_{np}$.
\begin{theoreme} Le Frobenius $F_*$ est un isomorphisme de la composante $p$-primaire de
$\H^i(\O_n^{\star})$ sur la composante $p$-primaire de
$p\H^i(\O_{pn}^{\star})$.
\end{theoreme}
\par\noindent
En effet, le diagramme ci-dessus est un morphisme de suites
exactes longues. La suite exacte des noyaux dit que la
multiplication par $p$ est un isomorphisme du noyau de $F_*$ sur
lui-même. Comme ce noyau est un groupe fini, le résultat en
découle.
\par\medskip
En fait, la connaissance de la suite spectrale de Bockstein permet
le calcul suivant~:
\begin{theoreme}
Soit $p$ un nombre premier. Le terme
$p^{k-1}\H^i(\O_n^{\star})/p^{k}\H^i(\O_n^{\star})$ du gradué de
la filtration $p$-adique de la cohomologie de de Rham de degré
total $n$ est naturellement isomorphe aux cocycles
dans $\O_{{n}/{p^k}}^{\ i}\otimes\Z/p\Z$,
si $i>0$ et $0<k\leq \nu_p(n)$, et il est nul sinon.
\end{theoreme}
\par\noindent
Ce r\'esultat se d\'emontre par une double r\'ecurrence descendante
sur $k$ et $i$, en utilisant les suites exactes~:
$$\xymatrix{
0\to p^{k-1}\H^i(\O_n^{\star})/p^{k}\H^i(\O_n^{\star})\to E^i_k
\ar[r]^-{\partial}&p^{k-1}\H^{i+1}(\O_n^{\star})
}$$
et l'identification de la page $E_k$ avec le complexe de de Rham
$$E_k\cong \O_{{n}/{p^k}}\otimes\Z/p\Z\ ,\ \ 0<k\leq \nu_p(n)\ .$$
\par\medskip
\begin{exemple}{\rm
En degr\'e total 4, la cohomologie de de Rham est nulle sauf en degr\'e
cohomologique 1 et 2. Le th\'eor\`eme pr\'ec\'edent donne~:
$\H^2(\O^4)\cong\La^2/2$ comme foncteur de $L$. Le terme $\H^1$ est
plus int\'eressant.
Modulo 2, il est \'egal au $\Z/2\Z$-vectoriel des cocycles dans $\O_2^1\otimes\Z/2\Z$, qui n'est
autre que la seconde puissance divis\'ee $\Gamma^2$ modulo 2. On a donc
une extension de foncteurs~:
$$0\to S^1/2\to\H^1(\O_4)\to\Gamma^2/2\to 0$$
qui fait appara\^\i tre $\H^1(\O_4)$ comme une variante des vecteurs de
Witt de longueur 2. En fait, cette extension repr\'esente la seule
classe d'extensions de foncteurs non scind\'ee entre $\Gamma^2/2$ et $S^1/2$
\cite[paragraphe 9.3]{FLS}.
}\end{exemple}
\bigskip\vfill
\bibliographystyle{amsalpha}

\eject
\selectlanguage{english}
\section*{Integral de {R}ham Cohomology}
\begin{abstract}
The Cartier isomorphism allows a nice description of the Bockstein
spectral sequence of the de Rham complex over the integers. It is
used to compute the integral de Rham cohomology of affine
spaces.\end{abstract}
\medskip
For a finite rank free abelian
group $L$, let $S^j(L)$ and $\La^i(L)$ respectively denote the
$j$-th symmetric and $i$-th exterior power of $L$. The de Rham
complex in total degree $n$ is the complex:
$$\O_n^{\star}\ :\ \xymatrix{S^n\ar[r]^-{d} &
  S^{n-1}\otimes\La^1\ar[r]^-{d}&\cdots\
  S^{n-i}\otimes\La^i\ar[r]^-{d}&\cdots\ar[r]^-{d}& \La^n\ar[r]^-{d}&
  0\ar[r]^-{d}&0 \ \cdots}$$
whose differential $d$ is the derivative of polynomials. The
purpose of this note is to describe its cohomology.
\begin{e-proposition}\label{e-torsion}
Every de Rham cohomology class in total degree $n$ is annihilated
by $n$.
\end{e-proposition}
\par\noindent
One way to prove that the de Rham cohomology of the affine space
over a characteristic zero field cancels is to use the Koszul
differential $\kappa$ as a contracting homotopy. Over the
integers, the formula $d\kappa +\kappa d=n$ only gives the above
proposition.
\par\medskip
The de Rham cohomology being a finite torsion group, we determine
it by computing the Bockstein spectral sequence of the de Rham
complex at each prime.
\begin{e-proposition}
Let $p$ be a prime. The Frobenius map $F_*$ induces a morphism
from the Bockstein exact couple to its derived couple. The map
$E_1\to E_2$ is an isomorphism, which is induced by the Cartier
isomorphism $C^{-1}$ :
$\O_n\otimes\Z/p\Z\to\H(\O_{pn}^{\star}\otimes\Z/p\Z)$. The map
$F_*$ : $\H^i(\O_n^{\star})\to p\H^i(\O_{pn}^{\star})$ yields an
isomorphism between the $p$-primary parts.
\end{e-proposition}
\par\medskip
It is amusing to note that the Bockstein spectral sequence is
stationary but it is not trivial:
\begin{theorem}
Each page of the Bockstein spectral sequence of the de Rham
complex is isomorphic to the de Rham complex over the prime field.
\end{theorem}
\begin{theorem}
Let $p$ be a prime. The de Rham cohomology in total degree $n$ has
a $p$-adic filtration whose factor
$p^{k-1}\H^i(\O_n^{\star})/p^{k}\H^i(\O_n^{\star})$ is naturally
isomorphic to the cocycles in $\O_{{n}/{p^k}}^{\ i}\otimes\Z/p\Z$,
if $i>0$ and $0<k\leq \nu_p(n)$, and else it vanishes.
\end{theorem}
\vfill
\end{document}